\documentclass[10pt]{article}
\usepackage{amsmath,amssymb,amsthm,amscd}
\numberwithin{equation}{section}

\newtheorem{prop}{Proposition}[section]
\newtheorem{theorem}[prop]{Theorem}

\newtheorem{remark}[prop]{Remark}
\newtheorem{example}[prop]{Example}
\newtheorem{definition}[prop]{Definition}
\newtheorem{conjecture}[prop]{Conjecture}
\newtheorem{acknowledgment}[prop]{Acknowledgment}

\def\<{\langle}
\def\>{\rangle}
\def\({\left(}
\def\){\right)}

\def\p{\partial}

\def\Ric{{\rm Ric}}

\begin{document}

\title{General Weak Limit for K\"ahler-Ricci Flow}

\author{the University of Sydney \\
Zhou Zhang  \footnote{Partially supported by Australian Research Council Discovery Project: DP110102654.}} 
\date{}
\maketitle

\begin{abstract}

Consider the K\"ahler-Ricci flow with finite time singularities over any closed K\"ahler manifold. We prove the existence of the flow limit in the sense of current towards the time of singularity. This answers affirmatively a problem raised by Tian in \cite{tian-survey} on the uniqueness of the weak limit from sequential convergence construction. The notion of minimal singularity in the study of positive current, introduced by Demailly, comes up naturally. We also provide some discussion on the infinite time singularity case for comparison. The consideration can be applied to more flexible evolution equation of K\"ahler-Ricci flow type for any cohomology class. The study is related to general conjectures on the singularities of K\"ahler-Ricci flows.

\end{abstract}

\section{Introduction}

Let $X$ be a closed K\"ahler manifold with $\dim_{\mathbb{C}}X=n\geqslant 2$. We consider the following K\"ahler-Ricci flow over $X$, 
\begin{equation}
\label{eq:krf}
\frac{\p\widetilde\omega (t)}{\p t}=-\Ric\(\widetilde\omega(t)\)-\widetilde\omega(t), ~~~~\widetilde\omega(0)=\omega_0. 
\end{equation}
where $\omega_0$ is the initial K\"ahler metric over $X$. This setting of K\"ahler-Ricci flow is general in the sense that we make no assumption on either the first Chern class of $X$, $c_1(X)$, or the initial K\"ahler class $[\omega_0]$ (not necessarily rational).  

The tremendous efforts and great successes in the study of this flow over the last decade are motivated by all kinds of intentions and visions from algebraic geometry, geometric analysis and several complex variables. The most prominent one would be Tian's Program, i.e., the Geometric Analytic Minimal Model Program. Simply speaking, the static equation for this flow (\ref{eq:krf}) is $\Ric\bigl(\widetilde\omega(\infty)\bigr)=-\widetilde\omega(\infty)$, which is in principle the desirable equation to study in search of the (singular and generalized) K\"ahler-Einstein with $K_X$ being almost positive (as in \cite{t-znote}, \cite{song-tian}, \cite{song-tian-2} and so on).   

\vspace{0.1in}

This general K\"ahler-Ricci flow can still be reduced to a scalar evolution equation for the metric potential as follows. Set $\omega_t=\omega_\infty+e^{-t}(\omega_0-\omega_\infty)$ with $\omega_\infty=-\Ric(\Omega)$ for some smooth volume form $\Omega$. The convention for $\Ric(\Omega)$ is that $\Ric(\omega^n)=\Ric(\omega)$. It's well known that $\widetilde\omega(t)=\omega_t+\sqrt{-1}\p\bar\p u$ with $u$ satisfying 
\begin{equation}
\label{eq:skrf}
\frac{\p u}{\p t}=\log\frac{(\omega_t+\sqrt{-1}\p\bar\p u)^n}{\Omega}-u, ~~~~u(\cdot, 0)=0,
\end{equation}
and this evolution equation is equivalent to (\ref{eq:krf}). 

By the optimal existence result in \cite{t-znote}, the smooth (metric) solution $\widetilde\omega_t$ exists as long as the class  
$$[\omega_t]=-c_1(X)+e^{-t}\bigl([\omega_0]+c_1(X)\bigr)=e^{-t}[\omega_0]+(1-e^{-t})\bigl(-c_1(X)\bigr)$$ 
remains in the K\"ahler cone of $X$, i.e., the open convex cone in the cohomology space $H^2(X; \mathbb{R})\cap H^{1,1}(X; \mathbb{C})$ consisting of all K\"ahler classes of $X$. Thus if we define {\it the time of singularity} 
$$T:=\sup\{t\,|\,[\omega_t]~\text{K\"ahler}\}\in (0, \infty],$$
the (classic) solution for the flow exists in $[0, T)$. For convenience, we apply the convention of $e^{-\infty}=0$. 

\vspace{0.1in}

If $K_X=-c_1(X)$ is K\"ahler, then $T=\infty$ and it's known that this flow always converges smoothly to the K\"ahler-Einstein metric for any $\omega_0$ (in any K\"ahler class). This is the non-degenerate case in \cite{t-znote} or more explicitly in \cite{thesis}. 

Otherwise, there has to be (metric) singularities developed for the flow metric $\widetilde\omega(t)$ as $t\to T$ because $[\omega_T]$ is on the boundary of the K\"ahler cone for $X$ (and no longer K\"ahler). The corresponding cases of $T<\infty$ and $T=\infty$ are naturally called the cases of finite and infinite time singularities respectively.  

Of course, the analysis of various metric singularities when approaching the time of singularity is crucial in understanding and dealing with the singularities for applications, for example, in order to carry out Tian's Program. For this purpose, it's useful to justify the existence and uniqueness of the limit, which provides the sole object for further discussion of regularity. For the uniqueness problem, the consideration in the weak sense is even more favorable.       

In general, as shown in Tian's survey \cite{tian-survey}, we always have a sequential weak limit for $\widetilde\omega(t)$ as $t\to T$ in the weak (i.e., current) sense. Let's point out that the argument in \cite{tian-survey} also works when $T=\infty$ for our version of the K\"ahler-Ricci flow because $\omega_t$ is uniformly controlled as form even if $T=\infty$. Normalization of the metric potential $u$ is performed in that work to achieve the sequential limit, and it's conjectured there that the sequential limit is actually unique, i.e., independent of the sequence chosen. 

This is the original motivation of this work. We attack this problem by studying the metric potential $u$ itself along the flow without normalization. If the flow limit exists, then the sequential limit is unique for sure.

In this work, we focus on the finite time singularity case, i.e., $T<\infty$, where we have the complete answer in the following theorem.

\begin{theorem}
\label{th:finite-time}
Over a closed K\"ahler manifold $X$, if the flow (\ref{eq:krf}) develops finite time singularities, i.e., the time of singularity $T<\infty$, then we have $\Phi=\omega_T+\sqrt{-1}\p\bar\p u(T)$ with $u(T)\in PSH_{\omega_T}(X)$ and $\widetilde\omega(t)\to\Phi$ in the weak sense as $t\to T$. In fact, for some positive constant $C$, $u+Ce^{-\frac{t}{2}}$ decreases to $u(T)+Ce^{-\frac{T}{2}}$ as $t\to T$, and $u(T)$ is of minimal singularities in $PSH_{\omega_T}(X)$. 
\end{theorem}

The situation for the infinite time singularity case is different in general as illustrated in Subsection 3.2. However, the above conclusion is true for the global volume non-collapsed case, even if $T=\infty$.  

\begin{prop}
\label{prop:non-collapsed}
Over a closed K\"ahler manifold $X$, consider the K\"ahler-Ricci flow (\ref{eq:krf}) with the time of singularity $T\in (0, \infty]$. If $[\omega_T]^n>0$, then we have $\Phi=\omega_T+\sqrt{-1}\p\bar\p u(T)$ with $u(T)\in PSH_{\omega_T}(X)$ and $\widetilde\omega(t)\to\widetilde\omega(T)$ in the weak sense as $t\to T$. In fact, for some positive constant $C$, $u+Ce^{-\frac{t}{2}}$ decreases to $u(T)+Ce^{-\frac{T}{2}}$ as $t\to T$, and $u(T)$ is of minimal singularities in $PSH_{\omega_T}(X)$. 
\end{prop}

The above conclusion when $T<\infty$ is weaker than the result in \cite{tosatti-collins}. When $T=\infty$, as shown in \cite{infinite-time-limit}, the inverse statement also holds. In this case, the metric potential actually stays uniformly bounded along the flow (and so is the limit $u(\infty)$) by the result in \cite{demailly-pali}, \cite{egz} and \cite{zzo}. The limit $u(\infty)$ is actually continuous by the result in \cite{zzo} and \cite{egz-viscosity}.   

\vspace{0.1in}

The organization of this work is as follows. In Section 2, we provide some general computations useful later and discussion some interesting special cases including Proposition \ref{prop:non-collapsed}. In Section 3, we prove Theorem \ref{th:finite-time}. In Section 4, we provide very brief but sufficient discussion on the generalization of Section 3 in the study of any cohomology class and elaborate on the relation with the notion of minimal singularity introduced by Demailly. Section 5 is about the relation between our consideration and general conjectures on the singularities for the K\"ahler-Ricci flow. 

\vspace{0.1in}

\noindent{\bf Notations:} as usual, we use $C$ to stand for a positive constant, possibly different at places. Also, $f\thicksim g$ means $\lim_{t\to T}\frac{f}{g}=1$.

\begin{acknowledgment}

The author would like to thank Professor G. Tian for introducing him to this wonderful field of research and constant encouragement. The interest of V. Tosatti in this work is very encouraging. It's invaluable for Professor S. Kolodziej to point out the regularization result by J. P. Demailly. The work was partially carried out during the visit to the School of Mathematical Sciences at Peking University and Beijing International Center for Mathematical Research, as part of the sabbatical leave from the University of Sydney. The author would like to thank Professor X. Zhu and these institutes for this wonderful opportunity.

\end{acknowledgment}

\section{The Plan and Special Cases}

In this section, we discuss the general plan and then summarize some interesting special cases. Although the argument in Section 3 takes care of them uniformly, they serve very well as motivations for this topic. Let's start with some general estimates for the K\"ahler-Ricci flow, leading to the plan of proving Theorem {\ref{th:finite-time}.     

The following computations are already quite standard as in \cite{t-znote} for example. The Laplacian, $\Delta$, is always with respect to the evolving metric along the flow, $\widetilde\omega(t)$.
\begin{equation}
\label{eq:t-derivative}
\frac{\partial}{\partial t}\(\frac{\partial u}{\partial t}\)=\Delta\(\frac{\partial u}{\partial t}\)-e^{-t}\<\widetilde\omega(t), \omega_0-\omega_\infty\>-\frac{\partial u}{\partial t},
\end{equation}
which is the $t$-derivative of (\ref{eq:skrf}). Here, $\<\omega, \alpha\>$ stands for the trace of the smooth real closed $(1, 1)$-form $\alpha$ with respect to the K\"ahler metric $\omega$. Equivalently, $\<\omega, \alpha\>=(\omega, \alpha)_\omega$ where $(\cdot, \cdot)_\omega$ is the hermitian inner product with respect to $\omega$. Then we transform (\ref{eq:t-derivative}) into the following two equations,
\begin{equation} 
\label{eq:1}
\frac{\partial}{\partial t}\(e^t\frac{\partial u}{\partial t}\)=\Delta\(e^t\frac{\partial u}{\partial t}\)-\<\widetilde\omega(t), \omega_0-\omega_\infty\>,
\end{equation}
\begin{equation}
\label{eq:2} 
\frac{\partial}{\partial t}\(\frac{\partial u}{\partial t}+u\)=\Delta\(\frac{\partial u}{\partial t}+u\)-n+\<\widetilde\omega(t), \omega_\infty\>.
\end{equation}
The difference of these two equations is
\begin{equation}
\label{eq:3}
\frac{\partial}{\partial t}\((1-e^t)\frac{\partial u}{\partial t}+u\)=\Delta\((1-e^t)\frac{\partial u}{\partial t}+u\)-n+\<\widetilde\omega(t), \omega_0\>.
\end{equation}
Since $\<\widetilde\omega(t), \omega_0\>>0$, by the standard Maximum Principle argument, it implies 
$$(1-e^t)\frac{\partial u}{\partial t}+u+nt\geqslant 0.$$
Also, there is a uniform upper bound for the metric potential $u$ by applying Maximum Principle to (\ref{eq:skrf}) directly. So this inequality can be rewritten as  
\begin{equation}
\label{ieq:u-t-derivative}
\frac{\partial u}{\partial t}\leqslant \frac{u+nt}{e^t-1}\leqslant \frac{nt+C}{e^t-1},
\end{equation}
which shows the essential decreasing (i.e., decreasing up to an exponentially decaying term) of $u$ along the flow. Notice that this control only depends on the upper bound of $u$. 

Taking $t$-derivative for (\ref{eq:t-derivative}), we arrive at 
\begin{equation}
\label{eq:t-derivatives}
\frac{\partial}{\partial t}\(\frac{\partial^2 u}{\partial t^2}\)=\Delta\(\frac{\partial^2 u}{\partial t^2}\)+e^{-t}\<\widetilde\omega(t),\omega_0-\omega_\infty\>-\frac{\partial^2 u}{\partial t^2}-{\vline\frac{\partial\widetilde\omega(t)}{\partial t}\vline}^2_{\,\widetilde\omega(t)}.
\end{equation}
Take the sum of (\ref{eq:t-derivative}) and (\ref{eq:t-derivatives}) to get
$$\frac{\partial}{\partial t}\(\frac{\partial^2 u}{\partial t^2}+\frac{\partial u}{\partial t}\)
=\Delta\(\frac{\partial^2 u}{\partial t^2}+\frac{\partial u}{\partial t}\)-\(\frac{\partial^2 u}
{\partial t^2}+\frac{\partial u}{\partial t}\)-{\vline\frac{\partial\widetilde\omega(t)}{\partial t}\vline}^2_{\,\widetilde\omega(t)}.$$
Applying Maximum Principle, we obtain  
$$\frac{\partial^2 u}{\partial t^2}+\frac{\partial u}{\partial t}\leqslant Ce^{-t},$$ 
which implies the essential decreasing of the volume form, $\widetilde\omega(t)^n=e^{\frac{\p u}{\p t}+u}\Omega$, along the flow, i.e.,  
$$\frac{\partial}{\partial t}\(\frac{\partial u}{\partial t}+u\)\leqslant Ce^{-t}.$$
From this inequality, we also have  
$$\frac{\partial u}{\partial t}\leqslant (C+Ct)e^{-t}\leqslant Ce^{-\frac{t}{2}}.$$
Thus, for some $C>0$, $u+Ce^{-\frac{t}{2}}$ is decreasing along the flow. By the well known property of plurisubharmonic function, as long as this term (or equivalently $u$) doesn't converge to $-\infty$ uniformly as $t\to T\in (0, \infty]$, $u$ will converge to some $u(T)\in PSH_{\omega_T}(X)$, and so $\Phi=\omega_T+\sqrt{-1}\p\bar\p u(T)$ is the flow limit for $\widetilde\omega(t)$ as $t\to T$ in the weak sense. Hence, in order to prove the existence of the flow limit in the weak sense for Theorem \ref{th:finite-time} and Proposition \ref{prop:non-collapsed}, we just need to {\it exclude the possibility of $u\to -\infty$ uniformly as $t\to T$}. This is our plan. In short, we don't normalize $u$ as in \cite{tian-survey}. 

\vspace{0.1in}

Next, we briefly recall the notion of minimal singularity introduced by Demailly in the study of positive $(1, 1)$-current.

\begin{definition}

Consider a smooth real closed $(1, 1)$-form $\alpha$ over $X$ with 
$$PSH_\alpha(X):=\{w\in L^1(X)\,|\, \alpha+\sqrt{-1}\p\bar\p w\geqslant 0~\text{in the sense of distribution}\}\neq\varnothing.$$ 
Then $u\in PSH_\alpha(X)$ is of minimal singularities if for any $v\in PSH_\alpha(X)$, $u\geqslant v-C$ for $C$ depending on $v$.  

\end{definition}

Any element in $PSH_\alpha(X)$ is well known to be bounded from above. So the singularity describes where and how the function approaches $-\infty$. Obviously, if $u\in PSH_\alpha(X)$ is actually bounded, then it's of minimal singularities.  

\begin{remark}
 
Since $|u|^p\widetilde\omega(t)^n=|u|^pe^{\frac{\p u}{\p t}+u}\Omega\leqslant C\Omega$ by the known upper bounds of $\frac{\p u}{\p t}$ and $u$, the weak convergence of all wedge powers $\widetilde\omega(t)^k\to\Phi^k$ as $t\to T$ is also available by the discussion in \cite{gz-07} as long as we have the existence of the weak limit $u(T)\in PSH_{\omega_T}(X)$ for the metric potential . 

\end{remark}

\subsection{Cases motivated by algebraic geometry interest}

In the consideration motivated by algebraic geometry interest, $[\omega_0]$ (or at least $[\omega_T]$) is rational, and the following claim is the key. 

\vspace{0.1in}

\noindent {\it Claim:} if $T<\infty$ and $[\omega_T]$ has a smooth real $(1, 1)$ non-negative form as representative, then the conclusion of Theorem \ref{th:finite-time} holds. 

\vspace{0.1in}

This is because we actually know $u\geqslant -C$ in this case by applying Maximum Principle argument directly on the flow. This is already done in \cite{t-znote} and we include some detail for later convenience. Start with  
\begin{equation}
\label{eq:4}
\frac{\partial}{\partial t}\((1-e^{t-T})\frac{\partial u}{\partial t}+u\)=\Delta\((1-e^{t-T})\frac{\partial u}{\partial t}+u\)-n+\<\widetilde\omega(t), \omega_T\>
\end{equation}
which is a proper linear combination of (\ref{eq:1}) and  (\ref{eq:2}). The assumption in the claim gives an $f\in C^\infty(X)$ such that $\omega_T+\sqrt{-1}\p\bar\p f\geqslant 0$. We then modify the above equation as follows. 
\begin{equation}
\begin{split}
\frac{\partial}{\partial t}\((1-e^{t-T})\frac{\partial u}{\partial t}+u-f\)
&= \Delta\((1-e^{t-T})\frac{\partial u}{\partial t}+u-f\) \\
&~~~~ -n+\<\widetilde\omega(t), \omega_T+\sqrt{-1}\p\bar\p f\>. \nonumber
\end{split}
\end{equation}
Applying Maximum Principle and noticing $T<\infty$, one has
$$(1-e^{t-T})\frac{\partial u}{\partial t}+u-f\geqslant -C.$$ 
Since $u\leqslant C$ and $\frac{\partial u}{\partial t}\leqslant C$, we conclude that for $t\in [0, T)$, 
$$u\geqslant -C, ~~~~\frac{\partial u}{\partial t}\geqslant -\frac{C}{1-e^{t-T}}\thicksim -\frac{C}
{T-t}.$$ 

In fact, one can get the lower bound for $u$ more directly using (\ref{eq:skrf}) as observed in \cite{song-weinkove}. After a proper choice of $\Omega$ when coming up with (\ref{eq:skrf}), we can make sure $\omega_T\geqslant 0$, and so $\omega_t\geqslant C(T-t)\omega_0$ for $t\in[0, T)$. Then by applying Maximum Principle to (\ref{eq:skrf}), we have
$$\frac{d \min_{X\times\{t\}}u}{d t}\geqslant n\log(T-t)-C-{\rm min}_{X\times\{t\}}u,$$
from which the lower bound of $u$ easily follows. 

\vspace{0.1in}

Together with the essential decreasing of $u$, we have the limit of the metric potential $u$ as a flow weak limit for $t\to T$, $u(T)\in PSH_{\omega_T}(X)\cap L^\infty(X)$. Thus in this case, by a well known result in pluripotential theory as in \cite{bed-tay}, one also has the flow weak convergence, $\widetilde\omega^k(t)\to\Phi^k=\bigl(\omega_T+\sqrt{-1}\p\bar\p u(T)\bigr)^k$ for $k=1, \cdots, n$. This $u(T)$ is bounded and so of minimal singularities. Hence Theorem \ref{th:finite-time} holds in this case.   

\vspace{0.1in}

Let's point out that if $[\omega_0]$ is a rational class (and so $X$ is projective for sure), the non-negativity of $[\omega_T]$ is available by the classic Rationality Theorem and the Base-Point-Free Theorem (for example in \cite{ko-mo} and \cite{kawa1}). So for algebraic geometry interest, the conclusion can already be drawn. Of course, our main purpose is to remove the assumption of $[\omega_0]$ being rational, and even $X$ being projective.    

\begin{remark}

It remains interesting to see whether the limit of $u$ is continuous, especially for the collapsed case, i.e.,  when $[\omega_T]^n=0$. For algebraic geometry interest, the continuity in the global volume noncollapsed case is known by \cite{zzo} and also \cite{egz-viscosity}. 

\end{remark}

The discussion here can be easily generalized to the case when $[\omega_T]-D$ is non-negative for an effective $\mathbb{R}$-divisor $D$. For simplicity of notations, we assume that $D$ is an effective $\mathbb{Z}$-divisor. Then $D$ can be identified as a holomorphic line bundle with a defining section $\sigma$ such that $D=\{\sigma=0\}$ and a hermitian metric $\|\cdot\|$. We get this information involved in the previous estimation as follows.
\begin{equation}
\begin{split}
\frac{\partial}{\partial t}\((1-e^{t-T})\frac{\partial u}{\partial t}+u-\log\|\sigma\|^2\)
&=\Delta\((1-e^{t-T})\frac{\partial u}{\partial t}+u-\log\|\sigma\|^2\)-n \\
&~~~~ +\<\widetilde\omega(t), \omega_T+\sqrt{-1}\p\bar\p\log\|\sigma\|^2\> \nonumber
\end{split}
\end{equation}
Since $[\omega_T]-D$ has a smooth non-negative representative, by choosing $\|\cdot\|$ properly, we have over $X\setminus D$,  
$$\omega_T+\sqrt{-1}\p\bar\p\log\|\sigma\|^2\geqslant 0.$$
So by Maximum Principle with the minimum clearly taken in $X\setminus D$ and noticing $T<\infty$, we arrive at 
$$(1-e^{t-T})\frac{\partial u}{\partial t}+u-\log\|\sigma\|^2\geqslant -C.$$ 
Again since $u\leqslant C$ and $\frac{\partial u }{\partial t}\leqslant C$, we conclude that 
$$u\geqslant \log\|\sigma\|^2-C, ~~\frac{\partial u}{\partial t}\geqslant\frac{C\log\|\sigma\|^2-C}{1-e^{t-T}}
\thicksim \frac{C\log\|\sigma\|^2-C}{T-t}.$$ 
This (degenerate) lower bound for $u$ is enough to guarantee the existence of its limit $u(T)\in PSH_{\omega_T}$ and the weak convergence of $\widetilde\omega(t)$ to $\Phi=\omega_T+\sqrt{-1}\p\bar\p u(T)$ as $t\to T$.

\begin{remark}

The case of $[\omega_T]=D$, an effective $\mathbb{R}$-divisor, is a special case in the above consideration. Also, this generalization would be more interesting in Section 4 when studying a general nef. (i.e., numerically effective) class using the more general evolution equation of K\"ahler-Ricci flow type.  

\end{remark}

\subsection{Global volume non-collapsed case}

Here, we show that if $[\omega_T]^n>0$, then $u$ can't go to $-\infty$ uniformly, and so $u\to u_T\in PSH_{\omega_T}(X)$. This provides the convergence statement in Proposition \ref{prop:non-collapsed}. Obviously, it is always the case that $[\omega_T]^n\geqslant 0$ since it is equal to the limit of $[\omega_t]^n>0$ as $t\to T$. Right now, we exclude the case of $[\omega_T]^n=0$, i.e., require the flow to be global volume non-collapsed. 

The proof is based on one simple observation. We rewrite (\ref{eq:skrf}) as follows. 
\begin{equation}
\label{eq:MA-type}
(\omega_t+\sqrt{-1}\p\bar\p u)^n=e^{\frac{\p u}{\p t}+u}\Omega.
\end{equation}
Assuming otherwise, $u$ converges to $-\infty$ uniformly over $X$ by the essential decreasing of $u$ and the basic property of plurisubharmonic function. 

Meanwhile, $[\omega_t]^n]\geqslant C>0$ for $t\in [0, T)$ since $[\omega_t]^n\to[\omega_T]^n>0$ as $t\to T$. We also know $\frac{\p u}{\p t}\leqslant C$. In light of 
$$\int_X e^u\Omega\geqslant C\int_X e^{\frac{\p u}{\p t}+u}\Omega=C\int_X(\omega_t+\sqrt{-1}\p\bar\p u)^n=C[\omega_t]^n\geqslant C>0,$$ 
we arrive at a contradiction. Notice that this argument works for $T\in (0, \infty]$. Hence, we conclude the weak convergence in this case. 

When $T=\infty$, then $[\omega_T]^n>0$ actually implies the uniform lower bound of $u$ by the result in \cite{demailly-pali}, \cite{egz} and \cite{zzo}. Clearly, the limit $u(T)$ is bounded and certainly of minimal singularities.

When $T<\infty$, if $[\omega_T]$ is rational (or for simplicity, $[\omega_0]$ is rational), it falls into the case considered in the previous subsection, so $u$ is bounded and of minimal singularities. When $[\omega_T]$ is a real (not necessary rational) class, we leave the justification of $u$ being of minimal singularities to the general discussion in Section 3.   

\begin{remark}

In the infinite time singularity case, as shown in \cite{infinite-time-limit}, the existence of unnormalized limit of $u$, i.e., $u$ not approaching $-\infty$ uniformly, indicates $K^n_X>0$, and so $u$ (and the limit) has to be uniformly bounded.  

\end{remark}

\begin{remark}
\label{remark-volume}
In general, one has $[\omega_t]^n\thicksim C(T-t)^k$ for some $k\in\{0, \cdots, n\}$ when $T<\infty$ and $[\omega_t]^n\thicksim Ce^{-kt}$ for some $k\in\{0, \cdots, n\}$ when $T=\infty$. This is seen by elementary consideration of the Taylor series of the explicit function $f(t)=[\omega_t]^n$ at $t=T$. In this subsection, $k=0$. All other $k$'s corresponds to the global volume collapsed case.  
\end{remark}

In the non-collapsed case, the following flow metric estimate, similar to that in \cite{scalar-r}, is of particular interest as described in Section 5. It begins with the inequality from parabolic Schwarz Lemma. Let $\phi=\<\widetilde\omega(t), \omega_0\>>0$. Using the computation for (\ref{eq:krf}) in \cite{song-tian}, one has 
\begin{equation}
\label{ieq:para-schwarz}
\(\frac{\p}{\p t}-\Delta\){\rm log}\phi\leqslant C_1\phi+1,
\end{equation}
where $C_1$ is a positive constant depending on the bisectional curvature of $\omega_0$. Also we have the following equation which is essentially (\ref{eq:3}): 
\begin{equation}
\label{eq:u-decreasing}
\(\frac{\p}{\p t}-\Delta\)\left((e^t-1)\frac{\p u}{\p t}-u-nt\right)=-\<\widetilde\omega(t), \omega_0\>.
\end{equation}
and so $(e^t-1)\frac{\p u}{\p t}-u-nt\leqslant 0$. Multiplying (\ref{eq:u-decreasing}) by a large enough constant $C_2>C_1+1$ and combining with (\ref{ieq:para-schwarz}), one arrives at 
\begin{equation}
\label{eq:main}
\begin{split}
\(\frac{\p}{\p t}-\Delta\)\left({\rm log}\phi+C_2\bigl((e^t-1)\frac{\p u}{\p t}-u-nt\bigr)\right)
&\leqslant nC_2+1-(C_2-C_1)\phi \\
&\leqslant C_3-\phi.
\end{split}
\end{equation}
Now apply Maximum Principle for the term ${\rm log}\phi+C_2\((e^t-1)\frac{\p u}{\p t}-u-nt\)$. Considering at the point where it achieves the maximal value, one has 
$$\phi\leqslant C,$$
and so 
$${\rm log}\phi+C_2\((e^t-1)\frac{\p u}{\p t}-u-nt\)\leqslant C,$$ 
which gives
$$\widetilde\omega(t)\leqslant Ce^{-C_2\((e^t-1)\frac{\p u}{\p t}-u-nt\)}\omega_0\leqslant Ce^
{-C(e^t\frac{\p u}{\p t}-t)}\omega_0.$$

Since $\widetilde\omega(t)^n=e^{\frac{\p u}{\p t}+u}\Omega$, we can further conclude that
$$Ce^{C(e^t\frac{\p u}{\p t}-t)}\omega_0\leqslant\widetilde\omega(t)\leqslant 
Ce^{-C(e^t\frac{\p u}{\p t}-t)}\omega_0.$$

If $T<\infty$, combining with the upper bound of $\frac{\p u}{\p t}$, we have for $t\in [0, T)$, 
\begin{equation}
\label{ieq:krf-metric-control}
Ce^{C\frac{\p u}{\p t}}\omega_0\leqslant\widetilde\omega(t)\leqslant Ce^{-C\frac{\p u}{\p t}}\omega_0.
\end{equation}
So the control of flow metric can be reduced to the lower bound of $\frac{\p u}{\p t}$. Although we know from \cite{scalar-r} that it's impossible to have a uniform lower bound for $\frac{\p u}{\p t}$, this control of metric is local and helpful in light of the local higher order estimates in \cite{sherman-weinkove}. This is useful for the discussion in Section 5. 

\begin{remark}

In the global collapsed case, this metric estimate blows up from both directions as $\frac{\p u}{\p t}\to\-\infty$ uniformly as $t\to T<\infty$, as discussed at the end of Subsection 3.2.   

\end{remark}

\subsection{Cases with curvature assumption}

Now we turn to several cases with various assumptions on Riemannian curvature tensor, all for the case of $T<\infty$. They help to motivate our consideration from Riemannian geometry point of view.   

\begin{itemize}

\item {\bf Case 1:} Ricci lower bound, i.e., $\Ric\big(\widetilde\omega(t)\bigr)\geqslant -C\widetilde\omega(t)$ for all $t\in [0, T)$. 

Clearly, $\widetilde\omega(t)\leqslant C\omega_0$ by the flow equation (\ref{eq:krf}). Thus we have 
$$-e^{-t}(\omega_0-\omega_\infty)+\sqrt{-1}\p\bar\p\frac{\p u}{\p t}=\frac{\p \widetilde
\omega(t)}{\p t}=-\Ric(\widetilde\omega(t))-\widetilde\omega(t)\leqslant C\widetilde\omega(t)\leqslant C\omega_0,$$
which gives
$$C\omega_0+\sqrt{-1}\p\bar\p\(-\frac{\p u}{\p t}\)\geqslant 0.$$

Applying the H\"ormander-Tian Inequality in \cite{tian-87}, there exist uniform constants $0<\alpha<1$ and $C>0$ such that
$$\int_X e^{\alpha\(\max_{_X}(-\frac{\p u}{\p t})-(-\frac{\p u}{\p t})\)}\Omega\leqslant C, ~\text{uniformly for any}~t\in[0, T),$$
which is    
$$\int_X e^{\alpha(-\min_{_X}\frac{\p u}{\p t}+\frac{\p u}{\p t})}\Omega\leqslant C.$$
So we have 
\begin{equation}
\label{ineq-3}
\int_X e^{\alpha\frac{\p u}{\p t}}\Omega\leqslant Ce^{\alpha\min_{_X}\(\frac{\p u}{\p t}\)}\leqslant Ce^{\alpha\frac{\p u}{\p t}\(x_{\min}(t), t\)}
\end{equation}
where $x_{\min}(t)$ is a point where $u(\cdot, t)$ takes the spatial minimum. Define the Lipschitz function $U(t):=\min_{X\times \{t\}} u$. We have $\frac{d U}{d t}=\frac{\p u}{\p t}(x_{\min}(t), t)$ in a proper sense.

Meanwhile, by Remark \ref{remark-volume}, we know 
$$\int_X e^{\frac{\p u}{\p t}+u}\Omega=[\omega_t]^n\geqslant C(T-t)^k$$
for some $k\in\{0, \cdots, n\}$. Together with $\alpha<1$ and the uniform upper bounds for $u$ and $\frac{\p u}{\p t}$, we arrive at  
\begin{equation}
\label{ineq-2}
\int_X e^{\alpha\frac{\p u}{\p t}}\Omega\geqslant C\int_X e^{\frac{\p u}{\p t}+u}\Omega\geqslant C(T-t)^k.
\end{equation}

Now combining (\ref{ineq-3}) with (\ref{ineq-2}), we have
$$\frac{d U}{d t}\geqslant C\log(T-t)-C$$
which gives $U\geqslant -C$. So there is a uniform $L^\infty$-bound for $u$. Let's summarize this in the following proposition. 

\begin{prop}
\label{prop:ricci-lower}
Consider the K\"ahler-Ricci flow (\ref{eq:krf}). If it develops finite time singularities with Ricci curvature uniformly bounded from below, then the metric potential in (\ref{eq:skrf}) has a uniform $L^\infty$-bound. When approaching the time of singularity, the flow metric weakly converges to a positive $(1, 1)$-current with bounded potential. 
\end{prop}

\begin{remark}

Using the same argument, one can replace the assumption of a uniform Ricci lower bound by $\Ric\bigl(\widetilde\omega(t)\bigr)\geqslant \alpha$ for a smooth $(1, 1)$-form $\alpha$, which is a priori a weaker assumption but less geometric. Meanwhile, by the result in \cite{ricci-lower}, we know that the assumption in the above proposition actually forces the global volume to collapse, i.e., $[\omega_T]^n=0$, which is indeed the difficult case in light of the discussion in Subsection 2.2.    

\end{remark}

\item {\bf Case 2:} Ricci form upper bound, i.e., $\Ric\bigl(\widetilde\omega(t)\bigr)\leqslant \alpha$ for a smooth $(1, 1)$-form $\alpha$.

This assumption is less geometric but still natural, for example in \cite{fong}. Notice that the upper bound $\Ric(\widetilde\omega(t))\leqslant C\widetilde\omega(t)$ indicates a positive lower bound for the flow metric for any finite time, and so actually rules out the finite time singularities as in \cite{scalar-r}. 

By the flow equation (\ref{eq:krf}), 
$$\alpha\geqslant \Ric\(\widetilde\omega(t)\)=-\frac{\p \widetilde\omega(t)}{\p t}-\widetilde\omega(t)=-\omega_\infty-\sqrt{-1}\p\bar\p\(\frac{\p u}{\p t}+u\).$$
So we have for some $C>0$, 
$$C\omega_0+\sqrt{-1}\p\bar\p\(\frac{\p u}{\p t}+u\)\geqslant 0.$$
The standard argument using Green's function then gives
$$\int_X\(\frac{\p u}{\p t}+u\)\omega^n_0\geqslant C\max\(\frac{\p u}{\p t}+u\)-C\geqslant C\log(T-t)-C$$
where the last step is easy to see since 
$$\int_Xe^{\frac{\p u}{\p t}+u}\Omega=[\widetilde\omega(t)]^n]=[\omega_t]^n\thicksim (T-t)^k$$
for some $k\in\{0, \cdots, n\}$. So we conclude 
$$\int_X u\omega^n_0\geqslant -C,$$
which prevents $u$ from going to $-\infty$ uniformly as $t\to T$. We leave the minimal singularity discussion to the general discussion in Section 3.  

\item {\bf Case 3:} Type I scalar curvature singularity , i.e., $|R(\widetilde\omega(t))|\leqslant\frac{C}{T-t}$. 

This is certainly weaker than the usual Type I singularity on Riemannian curvature. In light of the volume evolution equation 
$$\frac{\p \widetilde\omega^n_t}{\p t}=(-R-n)\widetilde\omega^n_t,$$
we have 
$$\frac{\p}{\p t}\(\frac{\p u}{\p t}+u\)\geqslant -\frac{C}{T-t}.$$
This gives $\frac{\p u}{\p t}+u\geqslant C\log(T-t)-C$, and so 
$$u\geqslant -C,$$
which is enough for the existence of the flow limit and convergence in the weak sense. The limit potential is bounded and certainly of minimal singularities. 

\end{itemize}

\section{The General Case}

Recall that in order to obtain the uniqueness of sequential limit, we consider instead the existence problem of flow weak limit. For this purpose, we need to rule out the possibility of $u\to -\infty$ uniformly over $X$ as $t\to T<\infty$. So far, we have justified this with various assumptions motivated by algebraic geometry and Riemannian geometry interests. Now we provide a general argument for the finite time singularity case. Then we illustrate the difference between the finite and infinite time singularity cases. 

\subsection{The proof of Theorem \ref{th:finite-time}}

We consider the case of $T<\infty$, i.e., the finite time singularity case. The following argument works in general. Of course, in some of the special cases discussed before, we have better estimates to draw the conclusion.   

To begin with, as $[\omega_T]$ is on the boundary of the K\"ahler cone for $X$, there exists $\varphi\in PSH_{\omega_T}(X)$. For example, one can use the sequential limit construction in \cite{tian-survey} to obtain such a function. 

Using the regularization result by Demailly (Theorem 1.6 in \cite{demailly-2014}), we have a decreasing approximation sequence for $\varphi$, $\{\varphi_m\}^\infty_{m=1}$, satisfying:

\begin{itemize}

\item $\varphi_m\in PSH_{\omega_T+\frac{1}{m}\omega_0}(X)$;

\item $\varphi_m\in C^\infty(X\setminus Z_m)$ with $Z_m\subset Z_{m+1}$ being analytic subvarieties of $X$, and along $Z_m$, $\varphi_m$ has logarithmic poles, i.e. locally being the logarithm of the sum of squares of holomorphic functions.

\end{itemize} 

We start with the following combination of (\ref{eq:3}) and (\ref{eq:4}),
\begin{equation} 
\begin{split}
&~~ \frac{\partial}{\partial t}\(\frac{1}{m}[(1-e^t)\frac{\partial u}{\partial t}+u]+[(1-e^{t-T})\frac{\partial u}{\partial t}+u]\) \\
&= \Delta\(\frac{1}{m}[(1-e^t)\frac{\partial u}{\partial t}+u]+[(1-e^{t-T})\frac{\partial u}{\partial t}+u]\) \\
&~~~~ -\frac{n(1+m)}{m}+\<\widetilde\omega(t), \frac{1}{m}\omega_0+\omega_T\>.
\end{split} \nonumber
\end{equation}

Then we modify it using $\varphi_m$: 
\begin{equation} 
\label{eq:5}
\begin{split}
&~~ \frac{\partial}{\partial t}\(\frac{1}{m}[(1-e^t)\frac{\partial u}{\partial t}+u]+[(1-e^{t-T})\frac{\partial u}{\partial t}+u]-\varphi_m\) \\
&= \Delta\(\frac{1}{m}[(1-e^t)\frac{\partial u}{\partial t}+u]+[(1-e^{t-T})\frac{\partial u}{\partial t}+u]-\varphi_m\) \\
&~~~~  -\frac{n(1+m)}{m}+\<\widetilde\omega(t), \frac{1}{m}\omega_0+\omega_T+
\sqrt{-1}\p\bar\p\varphi_m\>.
\end{split}
\end{equation}
where $ \frac{1}{m}\omega_0+\omega_T+\sqrt{-1}\p\bar\p\varphi_m$ is smooth and positive over $X\setminus Z_m$. As $\varphi_m\in C^\infty(X\setminus Z_m)$ and has $-\infty$ poles along $Z_m$, the spatial minimum over $X$ of 
$$\frac{1}{m}[(1-e^t)\frac{\partial u}{\partial t}+u]+[(1-e^{t-T})\frac{\partial u}{\partial t}+u]-\varphi_m$$ is always achieved in $X\setminus Z_m$, where everything is smooth. Now we can apply the standard Maximum Principle argument to conclude 
$$\frac{1}{m}[(1-e^t)\frac{\partial u}{\partial t}+u]+[(1-e^{t-T})\frac{\partial u}{\partial t}+u]-\varphi_m\geqslant -C,$$
which is uniformly for all $m$'s over $X\times [0, T)$. Here, the uniform upper bound of $\varphi_m$'s is used.  

We reformulate this inequality as follows:
$$\(1+\frac{1}{m}\)u+\(\frac{1}{m}(1-e^t)+(1-e^{t-T})\)\frac{\p u}{\p t}\geqslant -C+\varphi_m$$
For any $t\in[0, T)$, we could choose $m(t)$ large enough with 
$$0<\frac{2}{m(t)}(1-e^t)+(1-e^{t-T})<1,$$
where this choice is independent of other things. Combining with $u\leqslant C$, $\frac{\p u}{\p t}\leqslant C$ and $\varphi_{m(t)}\geqslant\varphi$, we arrive at 
$$u\geqslant -C+\varphi,$$
which is enough to exclude the possibility of $u\to -\infty$ uniformly as $t\to T$. 

This argument is for any $\varphi\in PSH_{\omega_T}(X)$, and so the flow limit $u(T)\in PSH_{\omega_T}(X)$ is of minimal singularities among all elements in $PSH_{\omega_T}(X)$. 

Hence we have completed the proof of Theorem \ref{th:finite-time}.

\begin{remark}

By Theorem \ref{th:finite-time}, the flow provides a smooth descreasing approximation for the finite time flow limit which is of minimal singularities. This interpretation will be more interesting in light of the discussion in Section 4. Also, there is similar discussion in the final version of \cite{tosatti-collins}, where the focus is certainly different. 

\end{remark}

\subsection{The difference between $T<\infty$ and $T=\infty$ cases}

The following example shows that for the infinite time (collapsed) case, i.e., $T=\infty$ and $[\omega_T]^n=0$, it is possible that $u\to-\infty$ uniformly as $t\to -\infty$.

\begin{example}

Suppose that $K_X=[\omega_\infty]$ gives a fibration structure of $X$ with the generic fibre of dimension $0<k\leqslant n$, i.e., $P: X\to\mathbb{CP}^N$ with $mK_X=P^*[\omega_{_{FS}}]$ and $P(X)$ of complex dimension $n-k$. Then the flow (\ref{eq:krf}) exists forever and $u\thicksim -kt$ as $t\to\infty$. 

This can be seen as follows. Begin with the following scalar potential flow
$$\frac{\p v}{\p t}=\log\frac{(\omega_t+\sqrt{-1}\p\bar\p v)^n}{\Omega}-v+kt, ~~~~v(\cdot, 0)=0.$$
Clearly, it still corresponds to the same metric flow (\ref{eq:krf}) and the relation between $u$ in (\ref{eq:skrf}) and $v$ is 
$$u=v+f(t)~~\text{with}~~\frac{d f}{d t}+f=-kt, ~~f(0)=0.$$
It's easy to get $f(t)\thicksim -kt$ and $\frac{d f(t)}{d t}\thicksim -k$ as $t\to\infty$. Rewrite the equation 
of $v$ as follows
$$(\omega_t+\sqrt{-1}\p\bar\p v)^n=e^{-kt}e^{\frac{\p v}{\p t}+v}\Omega$$
and then apply the $L^\infty$ estimates in \cite{egz2} and \cite{demailly-pali}, we have $|v|\leqslant 
C$ for all time. Hence $u\thicksim -kt$ which tends to $-\infty$ as $t\to\infty$.

By the result in \cite{song-tian-scalar}, we know $\frac{\p v}{\p t}$ and also $\frac{\p u}{\p t}$ are uniformly bounded. 

\end{example}

Indeed, assuming Abundance Conjecture, the assumption of the example always holds for $K_X$ nef. and $K^n_X=0$, and so the behavior of $u$ and $\frac{\p u}{\p t}$ is quite universal. Thus, the $T=\infty$ collapsed case needs to be treated differently. 

The difference can be understood in a very intuitive way. For the finite time collapsed case, one also has the following flow for $v$,
$$(\omega_t+\sqrt{-1}\p\bar\p v)^n=(T-t)^ke^{\frac{\p v}{\p t}+v}\Omega,$$
which corresponds to the same metric flow (\ref{eq:krf}) and 
$$u=v+f(t)~~\text{with}~~\frac{d f}{d t}+f=k\log(T-t), ~~f(0)=0.$$
Now we have $|f(t)|\leqslant C$ and $\frac{d f}{d t}\thicksim k\log(T-t)$. In principle, we expect that $v$ stays bounded or at least doesn't tend to $-\infty$ uniformly, and so that's also expected for $u$.

Meanwhile, the difference also exists regarding the behaviour of $\frac{\p u}{\p t}$. In fact, for the finite time collapsed case, we can justify $\frac{\p u}{\p t}\to -\infty$ uniformly as $t\to\ T$ by the following simple argument. By the discussion in Subsection 2.2, the limit of $\frac{\p u}{\p t}+u$ as $t\to T$ is $-\infty$ almost everywhere from the volume consideration. This limit is clearly essentially upper semi-continuous as the (essentially) decreasing limit of the smooth function $\frac{\p u}{\p t}+u$. So the limit of $\frac{\p u}{\p t}+u$ is indeed $-\infty$ over $X$, and the convergence is then uniform by the classic considersation as for Dini's theorem. Finally, by (\ref{ieq:u-t-derivative}), we conclude $\frac{\p u}{\p t}$ converges to $-\infty$ uniformly as desired. 

\section{Applications to Any Nef. Class}

In this section, we generalize the discussion in Section 3 to the study of a general nef. (i.e., numerically effective) class $\alpha$, which is a real $(1, 1)$-class on the boundary of the K\"ahler cone of $X$.     

For any K\"ahler metric $\omega_0$, we can choose a class $\beta\in H^{1, 1}(X; \mathbb{R})$ in the complement of the closure of the K\"ahler cone of $X$ such that the interval joining $[\omega_0]$ and $\beta$ intersects the boundary of the K\"ahler cone right at $\alpha$. Of course, the choice of $\beta$ is not unique even after fixing $[\omega_0]$ and $\alpha$, which is however not a concern for our purpose here. 

We then pick a smooth real $(1, 1)$-form $L$ representing $\beta$, and consider the following evolution equation of K\"ahler-Ricci flow type:
\begin{equation}
\label{eq:general-flow}
\frac{\p\widetilde\omega (t)}{\p t}=-\Ric\(\widetilde\omega(t)\)-\widetilde\omega(t)+\Ric(\Omega)+L, ~~~~\widetilde\omega(0)=\omega_0, 
\end{equation}
where $\Omega$ is a smooth volume form over $X$. This equation was considered in \cite{tsu-2} and further studied in \cite{t-znote} and \cite{thesis}. It shares a lot of common features as (\ref{eq:krf}), especially when considering the parabolic complex Monge-Amp\`ere equation for the metric potential.   

By the ODE consideration in $H^{1, 1}(X, \mathbb{R})$, we know $[\widetilde\omega(t)]=\beta+e^{-t}([\omega_0]-\beta)$, and the corresponding optimal existence result in \cite{t-znote} tells that the flow metric solution exists as long as the class $\beta+e^{-t}([\omega_0]-\beta)$ remains K\"ahler. So by our construction, the flow must develop finite time singularities. 
  
Set $\omega_t=L+e^{-t}(\omega_0-L)$. It's also known that $\widetilde\omega(t)=\omega_t+\sqrt{-1}\p\bar\p u$ with $u$ satisfying 
\begin{equation}
\label{eq:scalar-general-flow}
\frac{\p u}{\p t}=\log\frac{(\omega_t+\sqrt{-1}\p\bar\p u)^n}{\Omega}-u, ~~~~u(\cdot, 0)=0, 
\end{equation}
which is equivalent to (\ref{eq:general-flow}). 

All the computations and estimations at the beginning of Section 2 can be carried through in exactly the same way for this flow as for (\ref{eq:krf}). We still denote the time of singularity by $T<\infty$ with $\alpha=[\omega_T]=\beta+e^{-T}([\omega_0]-\beta)$. The construction in \cite{tian-survey} still works for this flow and so we know for example, $PSH_{\omega_T}(X)\neq\varnothing$. Then the same discussion as in Section 3 justifies the existence of the flow weak limit which is of minimal singularities. Hence we conclude: {\it using the evolution equation of K\"ahler-Ricci flow type (\ref{eq:general-flow}), one can construct a positive $(1, 1)$-current representative of minimal singularities for any nef. class together with a smooth approximation of it.}    

In light of the famous example by Serre as described in \cite{begz-big} about a nef. and big integral class without any positive current representative of bounded local potential, we know that the flow limit from the above construction doesn't always have bounded local potential. 

\vspace{0.1in}

In the following, we derive a flow metric estimate similar to that in Subsection 2.2. Let's begin with the following inequality mostly quoted from \cite{song-tian},
\begin{equation}
\(\frac{\p}{\p t}-\Delta\)\log\<\omega_0, \widetilde\omega(t)\>\leqslant C\<\widetilde\omega(t), \omega_0\>+C+\frac{\<\omega_0, \Ric(\Omega)+L\>}{\<\omega_0, \widetilde\omega(t)\>}, 
\end{equation}
where the last term on the right hand side comes from the extra term in (\ref{eq:general-flow}). Then one has
$$\(\frac{\p}{\p t}-\Delta\)\log\<\omega_0, \widetilde\omega_t\>
\leqslant C\<\widetilde\omega(t), \omega_0\>+C+\frac{C}{\<\omega_0, \widetilde\omega(t)\>}.$$

Let's recall (\ref{eq:3}) which still holds for (\ref{eq:general-flow}), 
$$\frac{\partial}{\partial t}\((1-e^t)\frac{\partial u}{\partial t}+u\)=\Delta\((1-e^t)\frac{\partial u}{\partial t}+u\)-n+\<\widetilde\omega(t), \omega_0\>,$$
which implies
$$(1-e^t)\frac{\p u}{\p t}+u+nt\geqslant 0.$$

Combining the two expressions above, we arrive at
\begin{equation}
\begin{split}
&~~ \(\frac{\p}{\p t}-\Delta\)\(\log\<\omega_0, \widetilde\omega(t)\>-B\bigl((1-e^t)\frac{\p u}{\p t}+u+nt\bigr)\) \\
&\leqslant (C-B)\<\widetilde\omega(t), \omega_0\>+C+\frac{C}{\<\omega_0, \widetilde\omega(t)\>}, \nonumber
\end{split}
\end{equation}
for a positive constant $B$ fixed later. At the (local space-time) maximum value point of the quantity, 
$$\log\<\omega_0, \widetilde\omega(t)\>-B\bigl((1-e^t)\frac{\p u}{\p t}+u+nt\bigr),$$ 
if it's not at the initial time (otherwise trivial), one has
$$(C-B)\<\widetilde\omega(t), \omega_0\>+C+\frac{C}{\<\omega_0, \widetilde\omega(t)\>}\geqslant 0.$$
Recall the elementary algebraic inequality
$$\<\omega_0, \widetilde\omega(t)\>\cdot\<\widetilde\omega(t), \omega_0\>\geqslant n^2,$$ 
and so one has $\frac{1}{\<\omega_0,\widetilde\omega(t)\>}\leqslant
\frac{\<\widetilde\omega(t), \omega_0\>}{n^2}$. After choosing $B>C+\frac{1}{n^2}+1$, at the maximal value point, we have
$$\<\widetilde\omega(t), \omega_0\>\leqslant C.$$
Now we apply another elementary inequality
$$\<\omega_0, \widetilde\omega(t)\>\leqslant\<\widetilde\omega(t), \omega_0\>^{n-1}\cdot\frac{\widetilde\omega(t)^n}{\omega^n_0}.$$ 
Together with $\widetilde\omega(t)^n=e^ {\frac{\p u}{\p t}+u}\Omega\leqslant C\Omega$, we have 
$$\<\omega_0, \widetilde\omega(t)\>\leqslant C$$
at that point. Noticing the lower bound for $(1-e^t)\frac{\p u}{\p t}+u+nt$, we conclude
$$\log\<\omega_0, \widetilde\omega(t)\>-B\bigl((1-e^t)\frac{\p u}{\p t}+u+nt\bigr)\leqslant C.$$
This gives
$$\widetilde\omega(t)\leqslant Ce^{B\((1-e^t)\frac{\p u}{\p t}+u+nt\)}\omega_0\leqslant Ce^
{-C(e^t\frac{\p u}{\p t}-t)}\omega_0.$$
Since $\widetilde\omega(t)^n=e^{\frac{\p u}{\p t}+u}\Omega$, we can further conclude that
$$Ce^{C(e^t\frac{\p u}{\p t}-t)}\omega_0\leqslant\widetilde\omega(t)\leqslant 
Ce^{-C(e^t\frac{\p u}{\p t}-t)}\omega_0.$$

Now we restrict to the case of $T<\infty$. Combining with the earlier estimates, we have for $t\in [0, T)$, 
\begin{equation}
\label{ieq:gkrf-metric-control}
Ce^{C\frac{\p u}{\p t}}\omega_0\leqslant\widetilde\omega(t)\leqslant Ce^{-C\frac{\p u}{\p t}}\omega_0.
\end{equation}
which is (\ref{ieq:krf-metric-control}) in this more general setting. 

Hence just as for (\ref{eq:krf}) as discussed in \cite{scalar-r}, there can not be any uniform lower bound for $\frac{\p u}{\p t}$ (or the volume form $\widetilde\omega^n(t)=e^{\frac{\p u}{\p t}+u}\Omega$) for the finite time singularity case.  

\begin{remark}

Similar to the result in \cite{sherman-weinkove}, with the uniform control of the flow metric, higher order estimates should be available even in a local fashion.  

\end{remark}

In all, the discussion so far for (\ref{eq:krf}) can be naturally generalized to the more general flow (\ref{eq:general-flow}). 

\section{Further Remarks}

Finally, we discuss the implication of this lower order consideration of the metric potential in understanding the formation of singularities for the K\"ahler-Ricci flow. Let's begin with the following conjecture. 

\begin{conjecture}

For the flow (\ref{eq:krf}) with singularity at $T<\infty$, $u\geqslant -C$ for $t\in [0, T)$. 

\end{conjecture}

Notice that the conjecture is on the classic K\"ahler-Ricci flow, i.e., about the canonical class $K_X$. For the more general flow discussed in Section 4, the situation is known to be different by Serre's example.  The confidence mostly comes from the cases from algebraic geometry consideration, in which the conjecture is known to hold as discussed in Subsection 2.1. Also, there is the natural relation with the conjectures on the singularity type of the K\"ahler-Ricci flow in \cite{song-weinkove-notes}, as illustrated by the discussion in Subsection 2.3.  

\vspace{0.1in}

Furthermore, there is this fundamental problem regarding the singularities of the K\"ahler-Ricci flow: {\it are they always developed along analytic varieties?} A little discussion with Professor F. Campana brought this to my attention.

In the global volume collapsed case, this would be the case if we consider the singularities in the usual sense for the unnormalized flow metric because of the vanishing of the global volume. In this case, it's more meaningful to search for global control of geometric quantities and a more precise understanding of the essential singularities regarding the metric collapsing, for example, after excluding the effect of the regular collapsing by proper scaling. 

So far, quite some progress has been made for the infinite time singularity case as in, for example, \cite{song-tian}, \cite{song-tian-2}, \cite{fong-z} and \cite{tosatti-yuguang}. The remaining difficulty is to achieve global geometric control in the presence of singular fibres. While for the finite time singularity case, it's fairly open with the existing results making serious assumptions, for example, in \cite{s-s-w} and \cite{fong}. 

The main result of this paper provides us with a flow weak limit of the metric potential which is closely related to the limiting class $[\omega_T]$. We expect the $-\infty$ locus of this limit of minimal singularities to characterize the essential singularities of the flow.    

\vspace{0.1in}

Meanwhile, the situation for the global volume non-collapsed case is a lot different as one could naturally expect the singularities to develop along a subvariety of $X$. There are evidences for both finite and infinite time singularity cases, for example, already in \cite{t-znote}, where the estimates are degenerate along a subvariety. 

Recently, there has been great progress made by Collins-Tosatti in the fundamental work \cite{tosatti-collins} on the finite time non-collapsed singularities of the K\"ahler-Ricci flow. More precisely, among other things, it's proved there that the flow stays smooth out of the subvariety $E_{nK}([\omega_T])={\rm Null}([\omega_T])$. This is done essentially by obtaining the proper lower bound of $\frac{\p u}{\p t}$. At the same time, it is impossible for the flow to stay smooth around any point of that subvariety in the Riemannian sense (i.e., with the curvature staying uniformly bounded in some fixed neighbourhood) by simple cohomology type consideration in light of the definition of the set ${\rm Null}([\omega_T])$ (i.e., the union of vanishing subvarieties with respect to $[\omega_T]$). It's quite obvious that the discussion in their work can be adapted to the more general flow (\ref{eq:general-flow}) in Section 4. Of course, for the higher order estimates, one needs to accept the statement in Remark 4.1.  
 
We now provide the following point of view which is certainly related but slightly different, coming from a variation of the definition for flow singularities. Let's consider the more general evolution equation of K\"ahler-Ricci flow type (\ref{eq:general-flow}) with finite time singularities, i.e., $T<\infty$.   

At the first sight, the (pluripolar) set $\{u_T=-\infty\}$, looks like the natural candidate for the singular set of the flow. However, this is not true in light of many known cases, for example, $u_T$ might well be bounded (and so this set is empty) in the presence of singularities. Indeed, by the discussion in Subsection 2.2, we are led to investigate the set 
$$\widehat S:=\{x\in X~|~\frac{\p u}{\p t}\to -\infty ~\text{for some time sequence}\}.$$ 
By (\ref{ieq:gkrf-metric-control}), the flow metric itself is (pointwise) bounded in the complement of $\widehat S$. 

\vspace{0.1in}

Since $\frac{\p u}{\p t}\leqslant Cu+C$ (from $\frac{\p u}{\p t}\leqslant\frac{u+nt}{e^t-1}$) and $u\leqslant C$, we have 
\begin{equation}
\label{ieq:equivalence}
C\frac{\p u}{\p t}-C\leqslant\frac{\p u}{\p t}+u\leqslant \frac{\p u}{\p t}+C.
\end{equation}
By the essential decreasing of the volume form, we can define the limit of $\frac{\p u}{\p t}+u$ which is an (essentially) upper semi-continuous function $V$ over $X$ valued in $[-\infty, C)$ for $C<\infty$, satisfying 
$$\int_X e^V\Omega=[\omega_T]^n.$$

By ({\ref{ieq:equivalence}}), we have 
$$\widehat S=\{x\in X~|~\frac{\p u}{\p t}\to-\infty\}=\{x\in X~|~\frac{\p u}{\p t}+u\to-\infty\}=\{x\in X~|~ V=-\infty\}.$$ 

In the global volume collapsed case, we have $V\equiv -\infty$ as discussed in Subsection 3.2, and so $\widehat S=X$ which coincides with the usual understanding of the singular set. 

Now we consider the global volume non-collapsed case. In the algebraic geometry setting, we already have in \cite{t-znote} that this set is contained in a subvariety of $X$. Indeed, by the estimate of Collins-Tosatti in \cite{tosatti-collins}, we have $\widehat S\subset{\rm Null}([\omega_T])$. 

The set $\widehat S$ can be complicated. More precisely, we have  
$$\widehat S=\cap^{\infty}_{A=1}\cup_{s\in [0, T)}\{x\in X~|~\frac{\p u}{\p t}+u+Ce^{-s}< -A~\text{at}~(x, s)\}.$$
The decreasing of $\frac{\p u}{\p t}+u+Ce^{-t}$ tells us that the open set 
$$\{x\in X~|~\frac{\p u}{\p t}+u+Ce^{-s}< -A~\text{at}~(x, s)\}$$ 
is increasing as $s\to T$. Also, the result in \cite{scalar-r} implies the open set 
$$\cup_{s\in [0, T)}\{x\in X~|~\frac{\p u}{\p t}+u+Ce^{-s}< -A~\text{at}~(x, s)\}\neq\varnothing$$
for any $A$. However, we are not even sure whether $\{x\in X~|~V=-\infty\}\neq \varnothing$, being the intersection of a sequence of decreasing open sets. A priori, $V$ might not actually take the value $-\infty$, although it can's have a lower bound by the discussion in \cite{scalar-r} for (\ref{eq:krf}) and the natural generalization to the general flow (\ref{eq:general-flow}) as described in Section 4. Nevertheless, if we consider the the lower semi-continuization $V_*$, the set 
$$S:=\{V_*=-\infty\}$$ 
is a closed set in $X$ and certainly non-empty by the above discussion. 

In the complement of $S$, $V$ is locally bounded (by the semi-continuity), and so the flow metric is locally uniformly bounded by (\ref{ieq:gkrf-metric-control}) and also for higher order estimates by the result in \cite{sherman-weinkove}, at least for the classic flow (\ref{eq:krf}). 

Thus, it is reasonable to consider $\{V_*=-\infty\}$ as the singular set. 

Let's point out that the function $V$ can a priori be wild. For example, $\{V=-\infty\}=\varnothing$ and $\{V_*=-\infty\}=X$ might happen simultaneously. However, by the result of Collins-Tosatti, $\widehat S\subset S\subset {\rm Null}([\omega_T])$. To conclude the discussion in this direction, we make the following conjecture.

\begin{conjecture} 

Consider the flow (\ref{eq:general-flow}) of K\"ahler-Ricci flow type. If there are finite time singularities with non-collapsed global volume, then in the notations above,  
$$\widehat S=S={\rm Null}([\omega_T]).$$   
\end{conjecture}

This predicts a more precise and also elemetary description for the blow-up along the singular set of the flow.

{\it 
\noindent Zhou Zhang \\ 
\noindent Address: Carslaw Building F07, the School of Mathematics and Statistics \\
the University of Sydney, NSW 2006, Australia \\
\noindent Email: zhangou@maths.usyd.edu.au \\
\noindent Fax: + 61 2 9351 4534}

\end{document}